\documentclass{article}

\usepackage{amsfonts,amssymb,amsmath}
\usepackage{graphicx}

\newtheorem{example}{Example}
\newtheorem{theorem}{Theorem}

\author{K. Igudesman and G. Shabernev}

\title{Novel method of fractal approximation}
\date{}

\begin{document}
\maketitle
\begin{abstract}
We introduce new method of optimization for finding free parameters
of affine iterated function systems (IFS), which are used for
fractal approximation. We provide the comparison of effectiveness of
fractal and quadratic types of approximation, which
are based on a similar optimization scheme, on the various types of
data: polynomial function, DNA primary sequence, price graph and
graph of random walking.
\end{abstract}

\section{Introduction}
It is well known that approximation is a crucial method for making
complicated data easier to describe and operate. In many cases we
have to deal with irregular forms, which can't be approximate with
desired precision. Fractal approximation become a suitable tool for
that purpose. Ideas for interpolation and approximation with the
help of fractals appeared in works of M.~Barnsley
\cite{barnsley_fractals_everywhere} and was developed by
P.~Massopust \cite{Massop} and C.~Bandt and A.~Kravchenko
\cite{bandt_kravchenko}.

Today we can apply fractals to approximate such interesting and
interdisciplinary data as graphs of DNA primary sequences of
different species and interbeat heart intervals \cite{landscapes},
price waves and many others.

Section \ref{S:fract_interp} of this work is devoted to the
construction of fractal interpolation functions. Necessary condition
on free parameters $d_i$ of affine iterated function systems is
shown. One graphical example is given.

In section \ref{S:approximation} we give the common scheme of
approximation of general function $g\in L^2[a,b]$ and obtain the
equation for direct calculation of free parameters $d_i$.

In section \ref{S:examples} we illustrate the results on concrete
examples.

\section{Fractal Interpolation Functions} \label{S:fract_interp}
There are two methods for constructing fractal interpolation
functions. In 1986 M.~Barnsley \cite{barnsley_fractals_everywhere}
defined such functions, as attractors of some specific iterated
function systems. In this work we use common approach, which was
developed by P.~Massopust \cite{Massop}.

Let $[a,b]\subset \mathbb{R}$ be a nonempty interval, $1<N\in\mathbb{N}$ and
$\{(x_i,y_i)\in[a,b]\times\mathbb{R}\mid a=x_0<x_1<\cdots<x_{N-1}<x_N=b\}$ --- are points of interpolation.
For all $i=\overline{1,N}$ consider affine transformations of the plane
$$
A_i:\mathbb{R}^2\rightarrow\mathbb{R}^2,\quad
A_i
\left(
\begin{array}{c}
  x \\
  y \\
\end{array}
\right)
:=
\left(
\begin{array}{cc}
  a_i & 0 \\
  c_i & d_i \\
\end{array}
\right)
\left(
\begin{array}{c}
  x \\
  y \\
\end{array}
\right)
+
\left(
\begin{array}{c}
  e_i \\
  f_i \\
\end{array}
\right).
$$
We require following two conditions hold true for all $i$:
$$
A_i(x_0,y_0)=(x_{i-1},y_{i-1}),\quad A_i(x_N,y_N)=(x_{i},y_{i}).
$$
In this case
\begin{equation}\label{E:coeff}
\begin{array}{ll}
\displaystyle a_i=\frac{x_i-x_{i-1}}{b-a},&
\displaystyle c_i=\frac{y_i-y_{i-1}-d_i(y_N-y_0)}{b-a},\\
\displaystyle e_i=\frac{bx_{i-1}-ax_i}{b-a},& \displaystyle
f_i=\frac{by_{i-1}-ay_i-d_i(by_0-ay_N)}{b-a},
\end{array}
\end{equation}
there $\{d_i\}_{i=1}^N$ act like family of parameters. Notice, that
for all $i$  operator $A_i$ takes the line segment between
$(x_0,y_0)$ and $(x_N,y_N)$ to the line segment passes through
points of interpolation $(x_{i-1},y_{i-1})$ and $(x_i,y_i)$.

Let
$\mathcal{K}$ be a space of nonempty compact subsets $\mathbb{R}^2$
with Hausdorff metric. Define the Hutchinson operator
\cite{hutchinson_fractals}
$$
\Phi:\mathcal{K}\rightarrow\mathcal{K},\qquad \Phi(E)=\bigcup_{i=1}^NA_i(E).
$$
It is easily seen \cite{barnsley_fractals_everywhere}, that the Hutchinson operator $\Phi$ take a graph of any
continuous function on a segment $[a,b]$ to a graph of a continuous function on the same segment.
Thus, $\Phi$ can be treated as operator on the space of continuous functions $C[a,b]$.

For all $i=\overline{1,N}$ denote
\begin{equation}\label{E:fun}
\begin{array}{l}
u_i:[a,b]\rightarrow[x_{i-1},x_i],\quad u_i(x):=a_ix+e_i,\\
p_i:[a,b]\rightarrow\mathbb{R},\quad p_i(x):=c_ix+f_i.
\end{array}
\end{equation}
Massopust \cite{Massop} has shown, that $\Phi$ acts on $C[a,b]$ according to the rule
\begin{equation}\label{E:operator}
(\Phi g)(x)=\sum_{i=1}^N\left((p_i\circ u_i^{-1})(x)+d_i(g\circ u_i^{-1})(x)\right)\chi_{[x_{i-1},x_i]}(x).
\end{equation}
Moreover, if
$|d_i|<1$ for all $i=\overline{1,N}$, then operator $\Phi$ is contractive
on the Banach space $(C[a,b],\|\ \|_{\infty})$ with contractive
constant $d\leq\max\{|d_i|\mid i=\overline{1,N}\}$.
By the fixed-point theorem there exists unique function $g^\star\in C[a,b]$, such that
$\Phi g^\star=g^\star$ and for all $g\in C[a,b]$ we have
$$
\lim_{n\to\infty}\|\Phi^n(g)-g^\star\|_\infty=0.
$$
We will call $g^\star$ fractal interpolation function.
It is clear, that if $g\in C[a,b]$, $g(x_0)=y_0$ and $g(x_N)=y_N$, then $\Phi(g)$ passes through points of interpolation.
In this case we will call $\Phi^n(g)$ pre-fractal interpolation functions of order $n$.

\begin{example}
Picture shows fractal interpolation function,
which was constructed on points of interpolation $(0,0)$, $(0.5,0.5)$ è $(1,0)$
with parameters $d_1=d_2=0.5$.
\begin{figure}[h!]
\centering \includegraphics{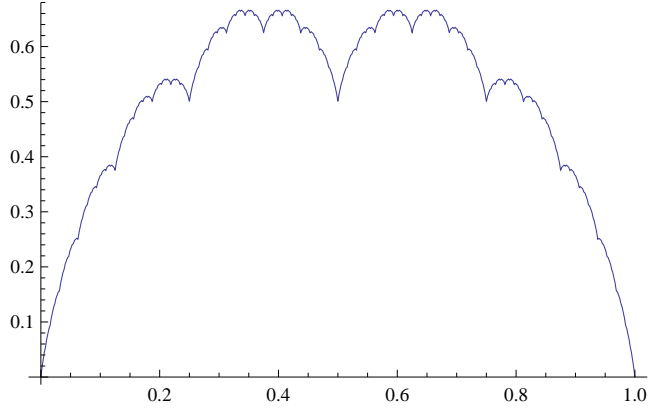}
\caption{Fractal interpolation function.} \label{F:ex1}
\end{figure}
\end{example}

\section{Approximation}\label{S:approximation}
From now on we assume, that $|d_i|<1$ for all $i=\overline{1,N}$. We
try to approximate function $g\in C[a,b]$ by the fractal
interpolation function $g^\star$, which is constructed on points of
interpolation $\{(x_i,y_i)\}_{i=0}^N$. Thus, it is sufficient to fit
parameters $d_i\in(-1,1)$ to minimize the distance between $g$ and
$g^\star$.

We use methods that have been developed for
fractal image compression \cite{barnsley_and_Hurd}.
Notice, that from (\ref{E:operator}), (\ref{E:fun}) and (\ref{E:coeff}) follows, that for all $g,h\in L^2[a,b]$
\begin{equation*}
\begin{split}
\|\Phi g-\Phi h\|_2=\sqrt{\int_a^b(\Phi g-\Phi h)^2\,\mathrm{d}x}
=\sqrt{\sum_{i=1}^Nd_i^2\int_{x_{i-1}}^{x_i}(g\circ u_i^{-1}(x)-h\circ u_i^{-1}(x))^2\,\mathrm{d}x}\\
\leq \max_{i=\overline{1,N}}\{|d_i|\}\cdot\sqrt{\sum_{i=1}^N a_i\int_a^b(g-h)^2\,\mathrm{d}x}
=\max_{i=\overline{1,N}}\{|d_i|\}\cdot\|g-h\|_2.
\end{split}
\end{equation*}
Thus, $\Phi:L^2[a,b]\rightarrow L^2[a,b]$ is contractive operator with a fixed point $g^\star$.

Furthermore, instead of minimization of $\|g-g^\star\|_2$ we will minimize $\|g-\Phi g\|_2$,
that makes the problem of optimization much easier. The collage theorem provides validity of such approach.
\begin{theorem}
Let $(X, d)$ be a non-empty complete metric space. Let
$T : X\to X$ be a contraction mapping on $X$ with contractivity factor $c<1$. Then for all $x\in X$
$$
d(x,x^\star)\leq\frac{d(x,T(x))}{1-c}
$$
where $x^\star$ is the fixed point of $T$.
\end{theorem}

$\blacktriangleright$
For all integer $n$ we have
\begin{equation*}
\begin{split}
d(x,x^\star)\leq d(x,T(x))+d(T(x),T^2(x))+\cdots +d(T^{n-1}(x),T^n(x))+d(T^n(x),x^\star)\\
\leq d(x,T(x))(1+c+c^2+\cdots+c^{n-1})+d(T^n(x),x^\star).
\end{split}
\end{equation*}
Letting $n\rightarrow\infty$ we establish the formula.
$\blacktriangleleft$

Considering (\ref{E:coeff}) and (\ref{E:fun}), we rewrite (\ref{E:operator}):
\begin{equation}\label{E:operator1}
(\Phi g)(x)=\sum_{i=1}^N\Big(\alpha_i(x)-d_i\big(\beta_i(x)-g\circ\gamma_i(x)\big)\Big)\chi_{[x_{i-1},x_i]}(x),
\end{equation}
where
$$
\alpha_i(x)=\frac{(y_i-y_{i-1})x+(x_iy_{i-1}-x_{i-1}y_i)}{x_i-x_{i-1}},$$
\begin{equation}\label{E:abg}
\beta_i(x)=\frac{(y_N-y_{0})x+(x_iy_{0}-x_{i-1}y_N)}{x_i-x_{i-1}},
\end{equation}
$$\gamma_i(x)=\frac{(b-a)x+(x_ia-x_{i-1}b)}{x_i-x_{i-1}}.$$
Thus, we have to minimize functional
$$
(\|g-\Phi g\|_2)^2=\sum_{i=1}^N\int_{x_{i-1}}^{x_i}\Big(g(x)-\alpha_i(x)+d_i\big(\beta_i(x)-g\circ\gamma_i(x)\big)\Big)^2\,\mathrm{d}x.
$$
Setting partial derivatives with respect to $d_i$ to zero we obtain
\begin{equation}\label{E:d}
d_i=\frac{\int_{x_{i-1}}^{x_i}\big(\alpha_i(x)-g(x)\big)\big(\beta_i(x)-g\circ\gamma_i(x)\big)\,\mathrm{d}x}
{\int_{x_{i-1}}^{x_i}\big(\beta_i(x)-g\circ\gamma_i(x)\big)^2\,\mathrm{d}x}, \qquad i=1,\ldots,N.
\end{equation}

\section{Discretization and results}\label{S:examples}
In this section we will approximate discrete data
$Z=\{(z_m,w_m)\}_{m=1}^M$, $a=z_1<z_2<\cdots <z_M=b$ by the fractal
interpolation function $g^\star$, which is constructed on points of
interpolation $X=\{(x_i,y_i)\}_{i=0}^N$, $N\ll M$. Taking $X\subset
Z$, $(x_0,y_0)=(z_1,w_1)$ and $(x_N,y_N)=(z_M,w_M)$ we fit
parameters $d_i\in(-1,1)$ to minimize
$$
\sum_{m=1}^M(w_m-g^\star(z_m))^2.
$$

Let us approximate $Z$ by the piecewise constant function $g:[a,b]\to\mathbb{R}$.
More precisely $g(z)=w_m$, where
$(z_m,w_m)\in Z$ and $z_m$ is a nearest neighbor of $z$.
From (\ref{E:d}) we obtain the discrete formulas for $d_i$:
\begin{equation}\label{E:d_discrete}
d_i=\frac{\sum\limits_{z_m\in[x_i,x_{i+1}]}\big(\alpha_i(z_m)-w_m\big)\big(\beta_i(z_m)-g\circ\gamma_i(z_m)\big)}
{\sum\limits_{z_m\in[x_i,x_{i+1}]}\big(\beta_i(z_m)-g\circ\gamma_i(z_m)\big)^2\,}, \qquad i=1,\ldots,N-1.
\end{equation}

After finding $d_i$ we obtain formulas for affine transformations
$A_i$ and we are able to construct fractal interpolation function
$g^\star$ for $g$.

Our aim is to compare fractal approximation with a piecewise
quadratic approximation function which is based on the same
discretization. On each segment $[x_{i-1},x_i]$ approximating
function has the quadratic form $q_i(x)=k_i x^2+r_i x+l_i$.
To get a continuous function we claim that $q_i(x_{i-1})=g(x_{i-1})$ and $q_{i}(x_{i})=g(x_{i})$.
From this we find coefficients $k_i$ and $l_i$. To find free parameter $r_i$ we
minimize functional
$$
\sum\limits_{z_m \in [x_{i-1},x_{i}]}(w_m-q_i(z_m))^2
$$
with respect to $r_i$ on each segment
$[x_{i-1},x_{i}],i=\overline{1,N}$. The approximating function
$q(x)$ will have following form:
$$
q(x)=\left\{
         \begin{array}{ll}
           q_1(x)=k_1 x^2+r_1 x+l_1, & x\in [a=x_0,x_1]; \\
           q_2(x)=k_2 x^2+r_2 x+l_2, & x\in [x_1,x_2]; \\
           \qquad \vdots\\
           q_N(x)=k_N x^2+r_N x+l_N, & x\in [x_{N-1},x_N=b].
         \end{array}
       \right.
$$
Since there is one free parameter $r_i$ in each function $q_i(x)$
and one parameter $d_i$ for each affine transformation $A_i$ it makes
the comparison correct.

To compare fractal and quadratic approximations we consider four types of data.
\begin{enumerate}
  \item Polynomial function.
  \item DNA sequence.
  \item Price graph.
  \item Random walking graph.
\end{enumerate}
For all types of data $M=10000$, $z_m=m$, $[a,b]=[1,M]$,
$\{w_m\}_{m=1}^M$ are normalized sequences, that is $E(\{w_m\})=0$
and $E(\{w_m^2\})=1$. For all cases we choose $(x_0,y_0)=(1,w_1)$,
$(x_N,y_N)=(M,w_M)$ and other interpolation points $(x_i,y_i)$,
$i=\overline{1,N-1}$ are local extremums of the given data.

\begin{example}
Let $f(x)=-6x + 5x^2 + 5x^3 - 5x^4 + x^5,\, x\in [-1,2.5]$. As we
work with the segment $[1,M]$ we map $[-1,2.5]$ to it. Consider
sequence $v_m=f\left(\frac{7(m-1)}{2(M-1)}-1\right)$,
$m=\overline{1,M}$. Set $w_m=(v_m-s_1)/s_2$, where $s_1$ and $s_2$
are mean and deviation of $\{v_m\}_{m=1}^M$. Figure
\ref{F:FuncGraph} shows the normalized sequence $\{w_m\}$.
\begin{figure}[h!]
\centering \includegraphics{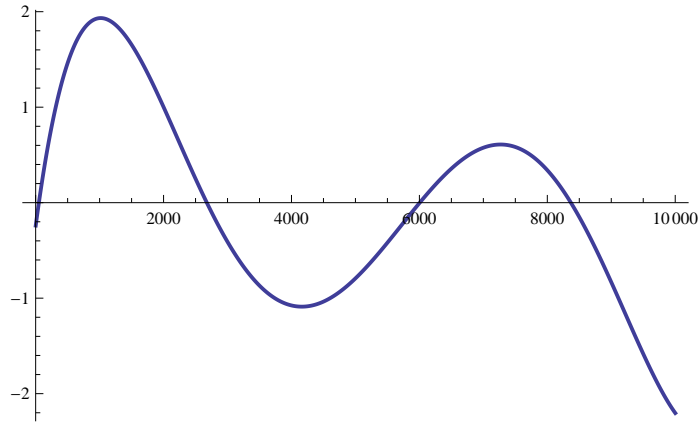}
\caption{The graph of original function $g$.} \label{F:FuncGraph}
\end{figure}
Choose five interpolation points $x_0=1$, $x_1=500$, $x_2=4000$,
$x_3=7500$, $x_4=10000$. Applying (\ref{E:d_discrete}) we obtain
$d_1=0.066$, $d_2=0.155$, $d_3=0.033$, $d_4=0.096$. The small values
of $|d_i|$ mean that on segments $[x_{i-1},x_i]$ fractal
approximation function looks as a straight line. Figure
\ref{F:FuncInterp} shows the graphs of fractal and quadratic
approximating functions.
\begin{figure}[h!]
\centering \includegraphics[width = 4.7in]{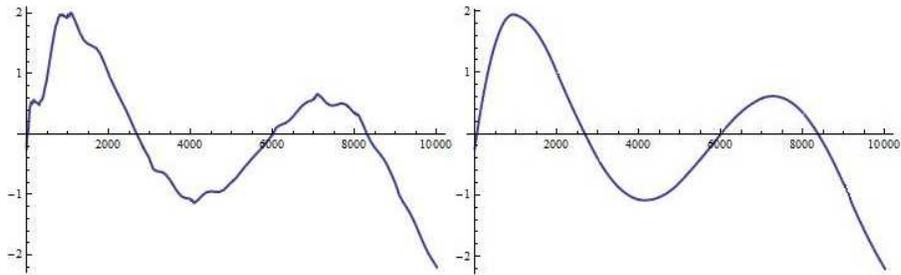}
\caption{Fractal and quadratic interpolations of the polynomial
function.}\label{F:FuncInterp}
\end{figure}
\end{example}

\begin{example}
A DNA sequence can be identified with a word over an alphabet
$\mathcal{N}=\{A,C,G,T\}$. Here we have the sequence of 10000
nucleotides of Edwardsiella tarda. The graph represented by the
formula $$v_1=0, \, v_m=v_{m-1}+\left\{
                             \begin{array}{ll}
                               +1, & \mbox{if} \,\, m^{th} \mbox{nucleotide belongs to (A,G)};\\
                               -1, & \mbox{if} \,\, m^{th} \mbox{nucleotide belongs to (C,T)}.
                             \end{array}
                           \right.$$
For full description of representation of DNA primary sequences see
\cite{prime_sequence}. Figure \ref{F:DNAGraph} shows the sequence
$\{w_m\}$ after normalization of $\{v_m\}$ according to the formula
in the previous example.
\begin{figure}[h!]
\centering \includegraphics{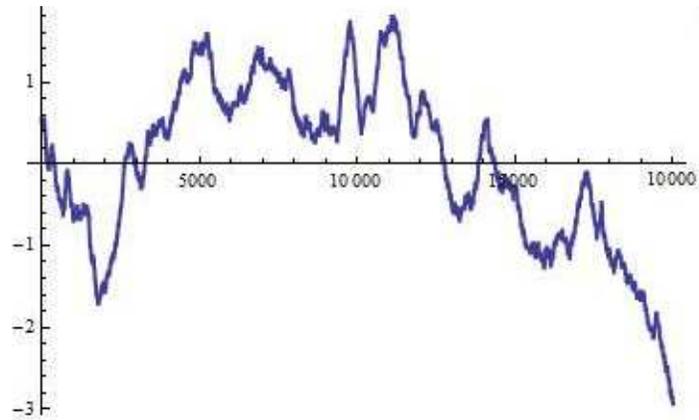} \caption{Picture shows DNA
Graph of 10000 nucleotides of Edwardsiella tarda.}
\label{F:DNAGraph}
\end{figure}
Interpolation points are $x_0=1$, $x_1=1000$, $x_2=2500$,
$x_3=3000$, $x_4=3500$,$x_5=5000$, $x_6=6500$, $x_7=7000$,
$x_8=8000$, $x_9=9000$, $x_{10}=10000$. Applying
(\ref{E:d_discrete}) we obtain $d_1=-0.001$, $d_2=0.274$,
$d_3=0.31$, $d_4=0.24$, $d_5=-0.057$, $d_6=0.211$, $d_7=-0.42$,
$d_8=-0.121$, $d_9=0.215$, $d_{10}=0.158$. Figure \ref{F:DNAInterp}
shows the graphs of fractal and quadratic approximating functions.
\begin{figure}[h!]
\centering \includegraphics[width = 4.7in]{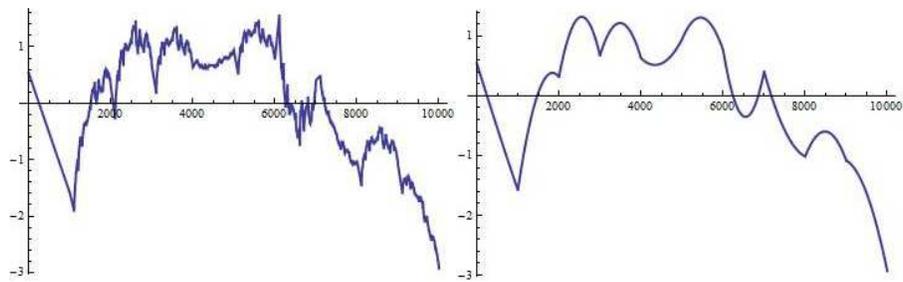}
\caption{Fractal and quadratic interpolations of the DNA
Graph.}\label{F:DNAInterp}
\end{figure}
\end{example}

\begin{example}
We take price wave of 10000 prices $v_m,\,m=\overline{1,M}$ of one
day period for EUR/USD, then normalize it
(Figure~\ref{F:PriceGraph}).
\begin{figure}[h!]
\centering \includegraphics{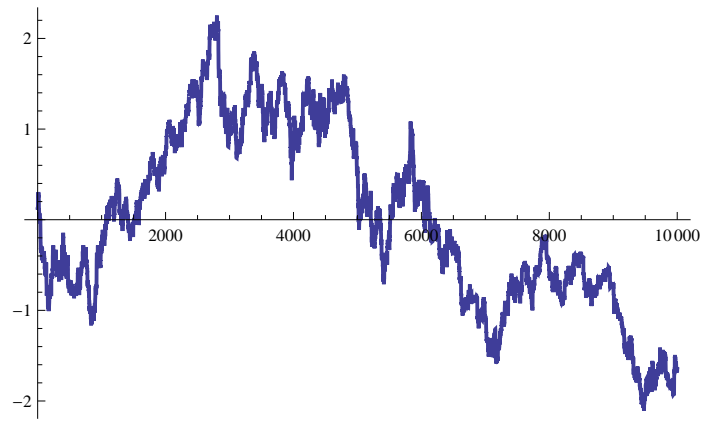}
\caption{Picture shows Price Graph for EUR/USD.} \label{F:PriceGraph}
\end{figure}
Interpolation points are $x_0=0$, $x_1=500$, $x_2=1500$, $x_3=2000$,
$x_4=2500$,$x_5=3000$,$x_6=4000$, $x_7=5000$, $x_8=6000$,
$x_9=8000$, $x_{10}=10000$. Applying (\ref{E:d_discrete}) we obtain
$d_1=-0.334$, $d_2=-0.004$, $d_3=0.315$, $d_4=0.307$, $d_5=0.333$,
$d_6=-0.28$, $d_7=-0.067$, $d_8=0.027$, $d_9=0.047$, $d_{10}=-0.33$.
Figure \ref{F:PriceInterp} shows the graphs of fractal and quadratic
approximating functions.
\begin{figure}[h!]
\centering \includegraphics[width = 4.7in]{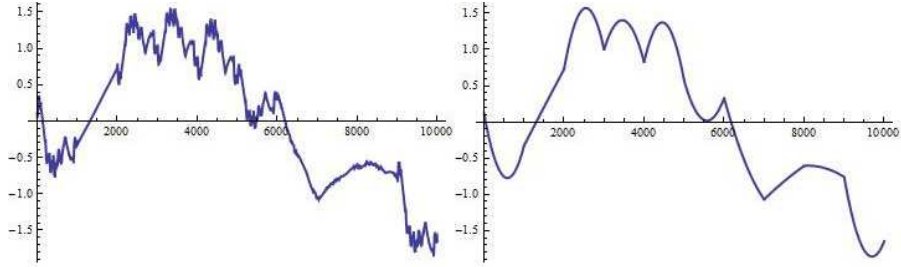}
\caption{Fractal and quadratic interpolations of the Price Graph.}
\label{F:PriceInterp}
\end{figure}
\end{example}

\begin{example}
Picture shows Random Walking Graph. It represented by the formula
$v_0=0, \, v_i=v_{i-1}+\xi_i$, where $\xi_i$ is a random value with
normal distribution.
\begin{figure}[h!]
\centering \includegraphics{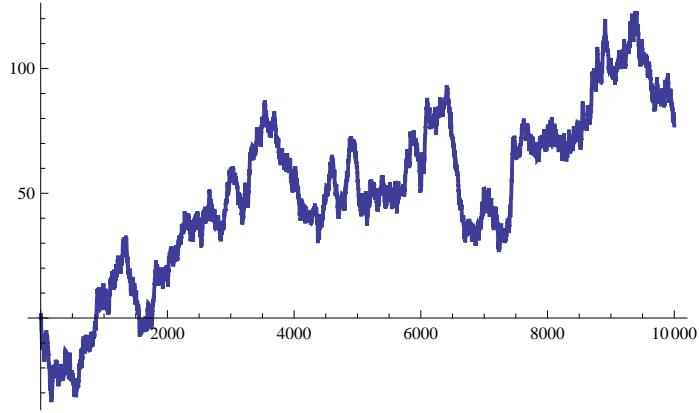}
\caption{Normalized Random Walking graph.} \label{F:RandomGraph}
\end{figure}
Interpolation points are $x_0=0$, $x_1=1500$, $x_2=2000$,
$x_3=3000$, $x_4=4000$, $x_5=5500$, $x_6=6300$, $x_7=7600$,
$x_8=8000$, $x_9=9000$, $x_{10}=10000$. Applying
(\ref{E:d_discrete}) we obtain $d_1=-0.237$, $d_2=0.14$,
$d_3=-0.020$, $d_4=-0.105$, $d_5=0.105$, $d_6=0.0545$, $d_7=-0.184$,
$d_8=-0.368$, $d_9=0.081$, $d_{10}=-0.111$. Figure
\ref{F:RandomWalkingInterp} shows the graphs of fractal and
quadratic approximating functions.
\begin{figure}[h!]
\centering \includegraphics[width = 4.7in]{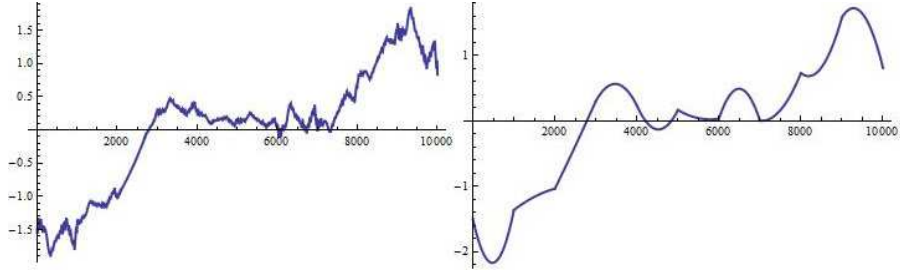}
\caption{Fractal and quadratic interpolations of the Random Walking
graph.} \label{F:RandomWalkingInterp}
\end{figure}
\end{example}

To compare the results we calculate approximation errors for each
type of data. Let $h(x)$ be the approximating function for data
$\{w_m\}_{m=1}^M$. Then approximation error is
$$
\sqrt{\sum\limits_{m=1}^M \frac{(h(x_m)-w_m))^2}{M}}.
$$

Here we represent the table of approximation errors for each type

\medskip
$$
\begin{array}{ccc}
           \quad & Fractal\,& Quadratic\,\\
  Polynomial Function & 0.0359037 & 0.0245094 \\
  DNA\, Primary\, Sequence & 0.0692072 & 0.0624714 \\
  Price\, Graph & 0.0501345 & 0.0533686  \\
  Random\, Walking & 0.1015339 & 0.101438
\end{array}
$$
\medskip

From it we see, that fractal approximation is better for price graph
and nearly equal for random walking, but much worse for smooth
function and slightly for DNA sequence. Different results were
appearing during calculations of errors. We assume that some
conditions could give us more exact approximation results from
fractal interpolation function and for that extra observations
should be established.

\end{document}